\documentclass{amsart}

\usepackage{amsmath,amsthm,amsfonts,amscd,amssymb}
\usepackage[pdfborder={0 0 0}]{hyperref}
\usepackage{verbatim} 

\newtheorem{thm}{Theorem}[section]

\newtheorem{lem}[thm]{Lemma}

\theoremstyle{definition}

\theoremstyle{remark}
\newtheorem{rem}{Remark}[section]

\begin{document}

\title[Lattice point counting and height bounds]{Lattice point counting and height bounds over number fields and quaternion algebras}
\author{Lenny Fukshansky and Glenn Henshaw}\thanks{The first author was partially supported by a grant from the Simons Foundation (\#208969 to Lenny Fukshansky) and by the NSA Young Investigator Grant \#1210223.}

\address{Department of Mathematics, 850 Columbia Avenue, Claremont McKenna College, Claremont, CA 91711}
\email{lenny@cmc.edu}

\address{Department of Mathematics, California State University at Channel Islands, One University Drive, Camarillo, CA 93012}
\email{glenn.henshaw@csuci.edu}


\begin{abstract}
An important problem in analytic and geometric combinatorics is estimating the number of lattice points in a compact convex set in a Euclidean space. Such estimates have numerous applications throughout mathematics. In this note, we exhibit applications of a particular estimate of this sort to several counting problems in number theory: counting integral points and units of bounded height over number fields, counting points of bounded height over positive definite quaternion algebras, and counting points of bounded height with a fixed support over global function fields. Our arguments use a collection of height comparison inequalities for heights over a number field and over a quaternion algebra. We also show how these inequalities can be used to obtain existence results for points of bounded height over a quaternion algebra, which constitute non-commutative analogues of variations of the classical Siegel's lemma and Cassels' theorem on small zeros of quadratic forms. 
\end{abstract}

\maketitle

\def\A{{\mathcal A}}
\def\AA{{\mathfrak A}}
\def\B{{\mathcal B}}
\def\C{{\mathcal C}}
\def\D{{\mathcal D}}
\def\E{{\mathcal E}}
\def\F{{\mathcal F}}
\def\x{{\mathcal H}}
\def\I{{\mathcal I}}
\def\J{{\mathcal J}}
\def\K{{\mathcal K}}
\def\kk{{\mathfrak K}}
\def\L{{\mathcal L}}
\def\M{{\mathcal M}}
\def\mm{{\mathfrak m}}
\def\MM{{\mathfrak M}}
\def\NN{{\mathfrak N}}
\def\N{{\mathcal N}}
\def\OO{{\mathfrak O}}
\def\O{{\mathcal O}}
\def\Pp{{\mathcal P}}
\def\pp{{\mathfrak p}}
\def\R{{\mathcal R}}
\def\s{{\mathcal S}}
\def\V{{\mathcal V}}
\def\UU{{\mathfrak U}}
\def\X{{\mathcal X}}
\def\Y{{\mathcal Y}}
\def\Z{{\mathcal Z}}
\def\H{{\mathcal H}}
\def\cee{{\mathbb C}}
\def\pee{{\mathbb P}}
\def\que{{\mathbb Q}}
\def\real{{\mathbb R}}
\def\zed{{\mathbb Z}}
\def\aaa{{\mathbb A}}
\def\hyp{{\mathbb H}}
\def\Bb{{\mathbb B}}
\def\Ff{{\mathbb F}}
\def\Nn{{\mathbb N}}
\def\kk{{\mathfrak K}}
\def\qbar{{\overline{\mathbb Q}}}
\def\kbar{{\overline{K}}}
\def\xbar{{\overline{x}}}
\def\ybar{{\overline{Y}}}
\def\kkbar{{\overline{\mathfrak K}}}
\def\ubar{{\overline{U}}}
\def\abar{{\overline{a}}}
\def\eps{{\varepsilon}}
\def\ahat{{\hat \alpha}}
\def\bhat{{\hat \beta}}
\def\gt{{\tilde \gamma}}
\def\h{{\tfrac12}}
\def\ba{{\boldsymbol a}}
\def\be{{\boldsymbol e}}
\def\bei{{\boldsymbol e_i}}
\def\bc{{\boldsymbol c}}
\def\bm{{\boldsymbol m}}
\def\bk{{\boldsymbol k}}
\def\bi{{\boldsymbol i}}
\def\bl{{\boldsymbol l}}
\def\bq{{\boldsymbol q}}
\def\bu{{\boldsymbol u}}
\def\bt{{\boldsymbol t}}
\def\bs{{\boldsymbol s}}
\def\bv{{\boldsymbol v}}
\def\bw{{\boldsymbol w}}
\def\bx{{\boldsymbol x}}
\def\bX{{\boldsymbol X}}
\def\bz{{\boldsymbol z}}
\def\bwy{{\boldsymbol y}}
\def\bY{{\boldsymbol Y}}
\def\bL{{\boldsymbol L}}
\def\bet{{\boldsymbol\eta}}
\def\bxi{{\boldsymbol\xi}}
\def\bo{{\boldsymbol 0}}
\def\bol{{\boldkey 1}_L}
\def\ep{\varepsilon}
\def\p{\boldsymbol\varphi}
\def\q{\boldsymbol\psi}
\def\Hf{H_{\fin}^{\O}}
\def\Hfo{H_{\fin}^{\O_1}}
\def\Hft{H_{\fin}^{\O_2}}
\def\Hfd{H_{\fin}^{O_D}}
\def\rank{\operatorname{rank}}
\def\aut{\operatorname{Aut}}
\def\lcm{\operatorname{lcm}}
\def\sgn{\operatorname{sgn}}
\def\spn{\operatorname{span}}
\def\md{\operatorname{mod}}
\def\Norm{\operatorname{Norm}}
\def\dim{\operatorname{dim}}
\def\det{\operatorname{det}}
\def\Vol{\operatorname{Vol}}
\def\rk{\operatorname{rk}}
\def\ord{\operatorname{ord}}
\def\ker{\operatorname{ker}}
\def\div{\operatorname{div}}
\def\Gal{\operatorname{Gal}}
\def\Tr{\operatorname{Tr}}
\def\nn{\operatorname{N}}
\def\inf{\operatorname{inf}}
\def\fin{\operatorname{fin}}
\def\Gr{\operatorname{Gr}}
\def\Mat{\operatorname{Mat}}
\def\GL{\operatorname{GL}}
\def\Ker{\operatorname{Ker}}
\def\Supp{\operatorname{Supp}}
\def\div{\operatorname{div}}
\def\Div{\operatorname{Div}}
\def\W{{\omega}}

\section{Introduction and statement of results}
\label{intro}

The classical combinatorial problem of estimating the number of lattice points in a compact set in the Euclidean space~$\real^N$, $N \geq 2$, has been studied extensively: see~\cite{lattice_points} for an overview of some of the main results. Estimates of this type have a great number of applications in many different areas of mathematics. In number theory and arithmetic geometry such results lead to the development of counting estimates for rational points on varieties over global fields.

A compact convex set in~$\real^N$ can be defined with the use of a norm, a device which measures ``size" of points in the space. Since~$\real$ is a local field, all norms on~$\real^N$ are equivalent. An analogous device over a global field is a height function, a standard tool of Diophantine geometry which measures size with respect to a full collection of infinitely many inequivalent norms simultaneously. A famous theorem of Northcott~\cite{northcott:ht} implies that any set of points of bounded height over a number field is finite. This observation is analogous to the statement that any compact set in~$\real^N$ contains only finitely many lattice points: in this more general case the number field plays the role of a lattice in the ambient adelic space, where inequalities on height define compact sets.

The first counting estimate on the number of algebraic numbers of bounded height in a fixed number field was produced by Schanuel~\cite{schanuel}. Schanuel's celebrated theorem has been extended and generalized in many ways by a number of authors over all global fields. While there are many further asymptotic results, extending Schanuel's original approach (see~\cite{masser_vaaler} and~\cite{widmer1} for some recent results and an overview), there are also several explicit bounds in the literature (see, for instance,~\cite{schmidt_northcott-1} and~\cite{loher_masser}). It should be remarked that only ~\cite{schmidt_northcott-1}  details some lower bounds, while the rest of the explicit estimates in the literature are upper bounds.

On the other hand, the problem of counting algebraic integers of bounded height in a fixed number field has received attention only more recently. While a mention of an asymptotic estimate without proof can be found in Lang's book~\cite{lang_dioph} (Theorem~5.2 on p.~70), to the best of our knowledge the first complete proofs of asymptotic estimates of this kind were obtained in~\cite{widmer2} and~\cite{barroero}. Explicit bounds in this situation are even more scarce, especially lower bounds. One explicit lower bound for the number of algebraic integers in a fixed number field was previously obtained by the first author in~\cite{null} (Corollary~1.6). Our first result is the following generalization of this bound; definition of the height function $h$ and other necessary notation is reviewed in Section~\ref{heights} below.

\begin{thm} \label{cnt_module} Let $K$ be a number field of degree $d$ over $\que$, $O_K$ its ring of integers, $N \geq 1$ an integer, and $\M \subset K^N$ a finitely generated $O_K$-module such that $\M \otimes_K K \cong K^L$, $1 \leq L \leq N$. Let $\D_K(\M)$ be the discriminant of the module $\M$, as given in~\eqref{module_disc} below. For a positive real number $R$, define
$$S_{K,N}(\M,R) = \left\{ \bx \in \M : h(\bx) \leq R \right\}.$$
Then
\begin{equation}
\label{mn_low_bnd}
\left|  S_{K,N}(\M,R) \right| \geq \left( \frac{R}{\E_1(K,\M,L) |\D_K(\M)|^{\frac{L}{2}}} - 1 \right) \left( \E_2(K,\M,L)R - 1 \right)^{L d-1},
\end{equation} 
for each 
$$R \geq \E_1(K,\M,L) |\D_K(\M)|^{L/2},$$
where constants $\E_1(K,\M,L)$ and $\E_2(K,\M,L)$ are defined in~\eqref{E1} and~\eqref{E2} below, respectively.
\end{thm}

\noindent
Our method of proof makes use of techniques in analytic and geometric combinatorics. Specifically, we employ the Minkowski embedding of the vector space $K^N$ into the Euclidean space~$\real^{Nd}$. The module~$\M$ under this embedding becomes a lattice of rank~$L d$, and the problem of counting points of bounded height in~$\M$ translates into the problem of counting lattice points in a certain compact domain in~$\real^{Nd}$. We then use a convenient explicit lattice point counting estimate in cubes as given by Lemma~\ref{cnt_lem} below. 

\begin{rem} \label{thm-1_rem} A simple upper bound on $|S_{K,N}(\M,R)|$ can be obtained from explicit estimates on the number of points of bounded height in $K^L$, as given in ~\cite{schmidt_northcott-1} and~\cite{loher_masser}.
\end{rem}

As an application of Theorem~\ref{cnt_module}, we obtain estimates on the number of points of bounded height which are integral over a fixed order in a positive definite quaternion algebra. To the best of our knowledge, this is the first application of lattice point counting techniques in a non-commutative situation. We start out by setting some basic notation. Let $K$ be a totally real number field of degree $d$ over $\que$, then $K$ has precisely $d$ real embeddings $\sigma_1,\dots,\sigma_d$. Let $O_K$ be the ring of integers in $K$ and let $\alpha,\beta \in O_K$ be totally negative elements, meaning that $\alpha^{(n)} := \sigma_n(\alpha) < 0$ and  $\beta^{(n)} := \sigma_n(\beta) < 0$ for all $1 \leq n \leq d$. Let $D = \binom{\alpha,\beta}{K}$ be a positive definite quaternion algebra over~$K$, generated by the elements $i,j,k$ which satisfy the following relations:
\begin{equation}
\label{quat_rel}
i^2=\alpha,\ j^2=\beta,\ ij=-ji=k,\ k^2=-\alpha\beta.
\end{equation}
It is possible to define height functions on $D$; we discuss definitions of three such heights in Section~\ref{heights} below: $h$, $H_{\inf}$, and $H^{\O}$, the last being a height function dependent on the choice of an order $\O$ is $D$. With this notation, we prove the following ``non-commutative analogue" of Theorem~\ref{cnt_module}.

\begin{thm} \label{main-1} Let $D = \binom{\alpha,\beta}{K}$ be as above and let $\O$ be an order in $D$. Let $N \geq 2$ be an integer, and let $Z \subseteq D^N$ be an $L$-dimensional right $D$-subspace, $1 \leq L \leq N$. For a positive real number $R$, define
$$S_{D,N}(Z,\O,R) = \left\{ \bx \in Z \cap \O^N : h(\bx) \leq R \right\}.$$
Then $\left| S_{D,N}(Z,\O,R) \right| \geq$
\begin{equation}
\label{count_ZO}
\left( \frac{R}{\E_3(D,\O,Z,d,L) H^{\O}(Z)^{4d}} - 1 \right) \left( \E_4(D,\O,Z,d,L) R - 1 \right)^{4L d-1},
\end{equation}
for each 
$$R \geq \E_3(D,\O,Z,d,L) H^{\O}(Z)^{4d},$$
where constants $\E_3(D,\O,Z,d,L)$ and $\E_4(D,\O,Z,d,L)$ are defined in~\eqref{E3} and~\eqref{E4} below, respectively.
\end{thm}
\smallskip

\noindent
To establish this result, we view $\O$ as an $O_K$-module, which allows us to apply Theorem~\ref{cnt_module}. Now the estimate is derived with the help of the height comparison lemmas proved in~\cite{quaternion}: these are inequalities relating heights over the number field $K$ to heights over the quaternion algebra $D$ over $K$. In fact, these inequalities can also be used to produce an upper bound on the number of points of bounded height in $D$ by an application of a result of~\cite{loher_masser}.

\begin{thm} \label{main-2} Let $D$ be as above, $R > 0$ be a real number, and define
\begin{equation}
\label{S_set}
S_{D,N}(R) = \{ \bx \in D^N : h(\bx) \leq R \}.
\end{equation}
Then
\begin{equation}
\label{count_bnd}
|S_{D,N}(R)| \leq (1088 d \log d)^{4N} \left( \frac{R}{t(\alpha,\beta)} \right)^{(4N+1)d},
\end{equation}
where the constant $t(\alpha,\beta)$ is defined below.
\end{thm}

\begin{rem} \label{north} Notice, in particular, that Theorem~\ref{main-2} implies Northcott's finiteness property for sets of points of bounded height on positive definite quaternion algebras over totally real number fields. Further, it is clear that $|S_{D,N}(R)| \geq |S_{D,N}(Z,\O,R)|$, which implies the upper bound of~\eqref{count_bnd} on $|S_{D,N}(Z,\O,R)|$. In addition,~\eqref{count_ZO} implies that
$$|S_{D,N}(R)| \gg_{N,K,D} R^{4Nd}.$$
\end{rem}
\smallskip

This paper is organized as follows. In Section~\ref{heights} we set the necessary notation, introduce height functions, and define the constants used in our estimates. We prove Theorems~\ref{cnt_module}, \ref{main-1}, and~\ref{main-2} in Section~\ref{count}. We also include two appendices with related results. In Appendix~\ref{S-unit} we show two more applications of the lattice point counting mechanism of Lemma~\ref{cnt_lem} to counting problems over global fields. Specifically, we obtain explicit estimates on the number of $S$-units of bounded height in an arbitrary number field as well as number of rational functions of bounded height supported on a given curve over a fixed finite field. Finally, in Appendix~\ref{points} we formulate a basic method (already used in deriving Theorem~\ref{main-1} from Theorem~\ref{cnt_module}) for obtaining results  over quaternion algebras by ``transferring" analogous results over number fields with the use of height comparison lemmas of~\cite{quaternion}. We further exhibit this method at work by obtaining existence results for points of bounded height in linear and quadratic spaces.
\bigskip

\section{Notation and heights}
\label{heights}

In this section we review the notation used in our main results and their proofs, as well as some further notation used in the appendices.

\subsection{Heights, quadratic forms, and constants over number fields.}
\label{ht_K}

Let $K$ be a number field of degree $d$ over $\que$, $O_K$ its ring of integers, $M(K)$ its set of places, $\D_K$ its discriminant, and let us write $\Nn$ for the norm from $K$ to $\que$. For each place $v \in M(K)$ we write $K_v$ for the completion of $K$ at $v$ and let $d_v = [K_v:\que_v]$ be the local degree of $K$ at $v$, so that for each $u \in M(\que)$
\begin{equation}
\sum_{v \in M(K), v|u} d_v = d.
\end{equation}

\noindent
For each place $v \in M(K)$ we define the absolute value $|\ |_v$ to be the unique absolute value on $K_v$ that extends either the usual absolute value on $\real$ or $\cee$ if $v | \infty$, or the usual $p$-adic absolute value on $\que_p$ if $v|p$, where $p$ is a rational prime. Then for each non-zero $a \in K$ the {\it product formula} reads
\begin{equation}
\label{product_formula}
\prod_{v \in M(K)} |a|^{d_v}_v = 1.
\end{equation} 

\noindent
We extend absolute values to vectors by defining the local heights. Let $N \geq 1$, and for each $v \in M(K)$ define a local height $H_v$ on $K_v^N$ by
$$H_v(\bx) = \max_{1 \leq i \leq N} |x_i|_v,$$
and for each $v | \infty$ define another local height $\H_v$ on $K_v^N$ by
$$\H_v(\bx) = \left( \sum_{i=1}^N |x_i|_v^2 \right)^{1/2}.$$
for each $\bx \in K_v^N$. Then we define two global height function on $K^N$:
$$H(\bx) = \prod_{v \in M(K)} H_v(\bx)^{d_v/d},\ \H(\bx) = \prod_{v \nmid \infty} H_v(\bx)^{d_v/d} \times \prod_{v | \infty} \H_v(\bx)^{d_v/d}$$
for each $\bx \in K^N$. Notice that due to the normalizing exponent $1/d$, our global height functions are absolute, i.e. for points over $\qbar$ their values do not depend on the field of definition. This means that if $\bx \in \qbar^N$ then $H(\bx)$ and $\H(\bx)$ can be evaluated over any number field containing the coordinates of $\bx$.

We also define an {\it inhomogeneous} height function on vectors by
$$h(\bx) = H(1,\bx),$$
hence $h(\bx) \geq H(\bx)$ for each $\bx \in \qbar^N$. In fact, the values of $H$ and $h$ are also related in the following sense: for each $\bx \in K^N$, there exists $a \in K$ such that $a\bx \in O_K^N$ and
\begin{equation}
\label{H_to_h}
H(\bx) = h(a \bx)
\end{equation}
when $N > 1$; when $N=1$, $h$ is just the usual Weil height.
\smallskip

We will also define two different height functions on matrices. First, let $B$ be an $N \times N$ matrix with entries in $K$, then we can view $B$ as a vector in $K^{N^2}$ and write $H(B)$ to denote the height of this vector. In particular, if $B$ is a symmetric matrix, then
$$Q(\bX,\bY) = \bX^t B \bY$$
is a symmetric bilinear form in $2N$ variables over $K$, and 
$$Q(\bX) := Q(\bX,\bX) =  \bX^t B \bX$$
is the associated quadratic form in $N$ variables. We define $H(Q)$, the height of such quadratic and bilinear forms, to be $H(B)$.
\smallskip

The second height we define on matrices is the same as height function on subspaces of $K^N$. Let $X = (\bx_1 \dots \bx_L)$ be an $N \times L$ matrix of rank $L$ over $K$, $1 \leq L \leq N$. Define
\begin{equation}
\label{matrix_ht}
\H(X) = \H(\bx_1 \wedge \dots \wedge \bx_L).
\end{equation}
For each $v | \infty$, the Cauchy-Binet formula guarantees that
\begin{equation}
\label{cauchy_binet}
\H_v(X) = |\det (X^* X)|_v^{1/2},
\end{equation}
where $X^*$ is the complex conjugate transpose of $X$. On the other hand, $\bx_1 \wedge \dots \wedge \bx_L$ can be identified with the vector $\Gr(X)$ of {\it Grassmann coordinates} of $X$ under the canonical embedding into $K^{\binom{N}{L}}$. Namely, let $\I$ be the collection of all subsets $I$ of $\{1,...,N\}$ of cardinality $L$, then $|\I| = \binom{N}{L}$. For each $I \in \I$, write $X_I$ for the $L \times L$ submatrix of $X$ consisting of all those rows of $X$ which are indexed by $I$. Define
\begin{equation}
\label{GR}
\Gr(X) = (\det (X_I))_{I \in \I} \in K^{\binom{N}{L}}.
\end{equation}
By our remark above, $\H(X) =  \H(\Gr(X))$. Now let $V \subseteq K^N$ be an $L$-dimensional subspace, $1 \leq L \leq N$. Choose a basis $\bx_1,...,\bx_L$ for $V$ over $K$, and let $X = (\bx_1\ ...\ \bx_L)$ be the corresponding $N \times L$ basis matrix. Define height of $V$ to be
$$H(V) := \H(X).$$
This height is well defined, since it does not depend on the choice of the basis for $V$: let $\bwy_1,...,\bwy_L$ be another basis for $V$ over $K$ and $Y = (\bwy_ 1 \dots \bwy_L)$ the corresponding $N \times L$ basis matrix, then there exists $C \in \GL_L(K)$ such that $Y = XC$, and so
$$\bwy_1 \wedge \dots \wedge \bwy_L = (\det C)\ \bx_1 \wedge \dots \wedge \bx_L,$$
hence, by the product formula $\H(\bwy_1 \wedge \dots \wedge \bwy_L) = \H(\bx_1 \wedge \dots \wedge \bx_L)$.

It will be convenient for us to define certain field constants that we use in our inequalities. First define 
\begin{equation}
\label{ckm}
c_K(\M) = \min \left\{ h(\alpha) : \alpha \in K \text{ such that } \alpha\M \subset O_K^L \right\},
\end{equation}
as well as
\begin{equation}
\label{zkm}
z_K(\M) = \min \left\{ h(\alpha) h(\alpha^{-1}) : \alpha \in K \text{ such that } \alpha\M \subset O_K^L \right\}.
\end{equation}
Now the constants used in the statement of Theorem~\ref{cnt_module} are given by
\begin{equation}
\label{E1}
\E_1(K,\M,L) = 2^{\frac{Lr_1-3}{2}} Ld\ z_K(\M) c_K(\M)^{Ld-1}
\end{equation}
and
\begin{equation}
\label{E2}
\E_2(K,\M,L) = \frac{2\sqrt{2}\ c_K(\M)}{Ld\ z_K(\M)}.
\end{equation}
Finally, for each $v | \infty$ and positive integer $j$ we define, as in \cite{vaaler:smallzeros},
\[ r_v(j) = \left\{ \begin{array}{ll}
    \pi^{-1/2} \Gamma(j/2+1)^{1/j} & \mbox{if $v | \infty$ is real,} \\
    (2\pi)^{-1/2} \Gamma(j+1)^{1/2j} & \mbox{if  $v | \infty$ is complex,}
\end{array}
\right. \]
and for any positive integers $\ell$ and $j$, define the constant $T_K(\ell, j)$ by
\begin{eqnarray}
\label{TK}
T_K(\ell,j) & = & 27 \left( \frac{1}{\pi} \right)^{\frac{r_2 \ell(9\ell+14)}{2d}} 2^{\frac{r_2\ell(9\ell+14) + (21\ell-21)d + 5r_1 + 4}{2d} + \max\{\ell, 9\}} \ell^{\frac{27\ell+51}{2}} j^{\frac{2}{d}} (j+2)^{\frac{3}{d}} \nonumber \\
& \times & |\D_K|^{\frac{\ell(9\ell+14) + 14}{2d} + \max\{\ell, 9\}} \left( \prod_{v | \infty} r_v(\ell-1)^{d_v/d} \right)^{\max\{\ell, 9\}}.
\end{eqnarray}
This constant is used in formula~\eqref{AK}, which is the definition of $\A_{K,\O} (L,M,J,\alpha,\beta)$, the constant in the inequality~\eqref{mn2} of Theorem \ref{main2}.

\subsection{Heights, quadratic forms, and constants over quaternion algebras.}
\label{ht_D}

We can also extend the height machinery to the context of quaternion algebras, using the approach of \cite{liebendorf:1}. Let $K$ be a totally real number field, $\alpha,\beta \in O_K$ be totally negative, and $D = \binom{\alpha,\beta}{K}$ be a positive definite quaternion algebra over~$K$, as defined in Section~\ref{intro} above. As a vector space, $D$ has dimension four over $K$, and $1,i,j,k$ is a basis. From now on we will fix this basis, and thus will always write each element $x \in D$ as
$$x = x(0)+x(1)i+x(2)j+x(3)k,$$
where $x(0),x(1),x(2),x(3) \in K$ are respective components of $x$, and the standard involution on $D$ is conjugation: 
$$\xbar= x(0)-x(1)i-x(2)j-x(3)k.$$
We define trace and norm on $D$ by
$$\Tr(x) = x+\xbar = 2x(0),\ \nn(x) = x\xbar = x(0)^2 - \alpha x(1)^2 - \beta x(2)^2 + \alpha \beta x(3)^2.$$
The algebra $D$ is said to be positive definite because the norm $\nn(x)$ is given by a positive definite quadratic form. In fact, since the norm form $\nn(x)$ is positive definite, $D_{v_n} := D \otimes_K K_{v_n}$ is isomorphic to the real quaternion $\hyp=\real + \real i + \real j + \real k$ for each $1 \leq n \leq d$. Hence each embedding $\sigma_n$ of $K$, $1 \leq n \leq d$, induces an embedding $\sigma_n : D \to D_{v_n}$, given by 
$$\sigma_n(x) = x(0)^{(n)} + x(1)^{(n)}i + x(2)^{(n)}j + x(3)^{(n)}k.$$
From now on we will write $x^{(n)}$ for $\sigma_n(x)$. Then the local norm at each archimedean place is also a positive definite quadratic form over the respective real completion~$K_{v_n}$:
\begin{eqnarray*}
& & \nn^{(n)}(x) = x^{(n)}\xbar^{(n)} \\
& & = \left(x(0)^{(n)}\right)^2 - \alpha^{(n)} \left(x(1)^{(n)}\right)^2 - \beta^{(n)} \left(x(2)^{(n)}\right)^2 + \alpha^{(n)} \beta^{(n)} \left(x(3)^{(n)}\right)^2,
\end{eqnarray*}
for each $1 \leq n \leq d$. We now have archimedean absolute values on $D$, corresponding to the infinite places $v_1,\dots,v_d$ of $K$: for each $x \in D$, define 
$$|x|_{v_n} = \sqrt{\nn^{(n)}(x)},$$
for every $1 \leq n \leq d$. It will be convenient to define
\begin{eqnarray}
\label{s_ab_v}
& & s_{v_n}(\alpha,\beta) = \max \{1,|\alpha|_{v_n},|\beta|_{v_n},|\alpha \beta|_{v_n}\}^{\frac{1}{2}}, \nonumber \\
& & t_{v_n}(\alpha,\beta) = \min \{1,|\alpha|_{v_n},|\beta|_{v_n},|\alpha \beta|_{v_n}\}^{\frac{1}{2}},
\end{eqnarray}
for each $1 \leq n \leq d$, and also let
\begin{equation}
\label{s_ab}
s(\alpha,\beta) = \prod_{n=1}^d s_{v_n}(\alpha,\beta),\ t(\alpha,\beta) = \prod_{n=1}^d t_{v_n}(\alpha,\beta).
\end{equation}
Since local norm forms are positive definite, we immediately have the following inequalities:
\begin{equation}
\label{loc_ineq}
t_{v_n}(\alpha,\beta) \max_{0 \leq m \leq 3} |x(m)|_{v_n} \leq\ |x|_{v_n} \leq 2 s_{v_n}(\alpha,\beta) \max_{0 \leq m \leq 3} |x(m)|_{v_n}.
\end{equation}
Now, generalizing notation of \cite{liebendorf:1}, we can define an infinite homogeneous height on $D^N$ by
\begin{equation}
\label{Hinf}
H_{\inf}(\bx) = \left( \prod_{n=1}^d \max_{1 \leq l \leq N} |x_l|_{v_n} \right)^{1/d},
\end{equation}
and define an infinite inhomogeneous height on $D^N$ by
\begin{equation}
\label{hinf}
h_{\inf}(\bx) = H_{\inf}(1,\bx),
\end{equation}
for every $\bx \in D^N$. Clearly, $H_{\inf}(\bx) \leq h_{\inf}(\bx)$. The infinite height takes into account the contributions at the archimedean places. As in \cite{liebendorf:1}, we also define its counterpart, the finite height. Let us once and for all fix an order $\O$ in $D$; our definition will be with respect to the order $\O$, and this height will be denoted by $\Hf$. Specifically, for each $\bx \in \O^N$, let
\begin{equation}
\label{Hfin_O}
\Hf(\bx) = [\O : \O x_1 + \dots + \O x_N]^{-1/4d}.
\end{equation}
This is well defined, since $\O x_1 + \dots + \O x_N$ is a left submodule of $\O$. Now we can define the global homogeneous height on $\O^N$ by
\begin{equation}
\label{HD}
H^{\O}(\bx) = H_{\inf}(\bx) \Hf(\bx),
\end{equation}
and the global inhomogeneous height by
\begin{equation}
\label{hD}
h(\bx) := H_{\inf}(1,\bx) \Hf(1,\bx) = h_{\inf}(\bx) \geq H^{\O}(\bx),
\end{equation}
since $\O + \O x_1 + \dots + \O x_N = \O$. To extend this definition to $D^N$, notice that for each $\bx \in D^N$ there exists $a \in O_K$ such that $a\bx \in \O^N$, and define $H^{\O}(\bx)$ to be $H^{\O}(a \bx)$ for any such~$a$. This is well defined by the product formula, and $H^{\O}(\bx t) = H^{\O}(\bx)$ for all $t \in D^{\times}$.
\smallskip

We will now define height on the set of proper right $D$-subspaces of $D^N$, again following \cite{liebendorf:1}. Recall that $D$ splits over $E=K(\sqrt{\alpha})$, meaning that there exists a $K$-algebra homomorphism $\rho: D \to \Mat_{22}(E)$, given~by
\begin{equation}
\label{rho}
\rho(x(0)+x(1)i+x(2)j+x(3)k) = \left( \begin{matrix} x(0)+x(1)\sqrt{\alpha} & x(2)+x(3)\sqrt{\alpha} \\ \beta(x(2)-x(3)\sqrt{\alpha}) & x(0)-x(1)\sqrt{\alpha} \end{matrix} \right),
\end{equation}
so that $\rho(D)$ spans $\Mat_{22}(E)$ as an $E$-vector space (see Proposition 13.2a (p. 238) and Exercise 1 (p. 240) of \cite{pierce}). This map extends naturally to matrices over $D$. Let $Z \subseteq D^N$ be an $L$-dimensional right vector subspace of $D^N$, $1 \leq L < N$. Then there exists an $(N-L) \times N$ matrix $C$ over $D$ with left row rank $N-L$ such that $Z$ is the solution space of the linear system $C\bX = \bo$.  Define
\begin{equation}
\label{Hinf_C}
H_{\inf}(C) = \left( \prod_{n=1}^d \left| \det\left( \rho(CC^*) \right) \right|_{v_n} \right)^{1/4d},
\end{equation}
where $C^*$ is the conjugate transpose of $C$. The analogue of Cauchy-Binet formula works here as well (see (2.7) and (2.8) of \cite{liebendorf:1}, as well as Corollary 1 of \cite{liebendorf:2}), and so we have an alternative formula:
\begin{equation}
\label{Hinf_C1}
H_{\inf}(C) = \left( \prod_{n=1}^d \sum_{C_0} \left| \det\left( \rho(C_0) \right) \right|^2_{v_n} \right)^{1/2d},
\end{equation}
where the sum is taken over all $(N-L) \times (N-L)$ minors $C_0$ of $C$. Also define
\begin{equation}
\label{Hfin_C}
\Hf(C) = [\O^{N-L} : C(\O^N)]^{-1/4d},
\end{equation}
where $C$ is viewed as a linear map $\O^N \to \O^{N-L}$. Then we can define
\begin{equation}
\label{H_Z}
H^{\O}(Z) = H^{\O}(C) := H_{\inf}(C) \Hf(C).
\end{equation}
This definition does not depend on the specific choice of such matrix~$C$. By the duality principle proved in \cite{liebendorf:3},
\begin{equation}
\label{H_Z_dual}
H^{\O}(Z) = H^{\O}(Z^{\perp}),
\end{equation}
where $Z^{\perp} = \{ \bwy \in D^N : \bx^* \bwy = 0\ \forall\ \bx \in Z \}$. This means that if $\bx_1,\dots,\bx_L$ is a basis for $Z$ over $D$ and $X = (\bx_1 \dots \bx_L)$ is the corresponding basis matrix, then 
\begin{equation}
\label{H_Z_dual1}
H^{\O}(Z) = H^{\O}(X) := \left( [\O^L : X^t(\O^N)]^{-1} \prod_{n=1}^d \left| \det\left( \rho(X^*X) \right) \right|_{v_n} \right)^{1/4d},
\end{equation}
completely analogous to the definition of the height $H^{\O}(C)$ in (\ref{H_Z}); here $X^t$ is viewed as a linear map $\O^N \to \O^{N-L}$.

It will also be convenient to define a map $[\ ]: D \to K^4$, given by 
$$[x] = (x(0),x(1),x(2),x(3)),$$
for each $x = x(0)+x(1)i+x(2)j+x(3)k \in D$. This map obviously extends to $[\ ]: D^N \to K^{4N}$, given by $[\bx] = ([x_1],\dots,[x_N])$ for each $\bx = (x_1,\dots,x_N) \in D^N$. Clearly this is a bijection; in fact, it is an isomorphism of $K$-vector spaces, and we will write $[\ ]^{-1}$ for its inverse. 
\smallskip

By analogy with heights over $D$, we will also write
$$H_{\inf}(\bx) = \prod_{v | \infty} H_v(\bx)^{d_v/d},\ H_{\fin}(\bx) = \prod_{v \nmid \infty} H_v(\bx)^{d_v/d},$$
for every $\bx \in K^N$. Then by Lemma 2.1 of \cite{liebendorf:1}, for every $\bx \in O_K^N$ we have
\begin{equation}
\label{hqf1}
H_{\fin}(\bx) = \left[ O_K : O_K x_1 + \dots + O_K x_N \right]^{-1/d}.
\end{equation}
Also, if $V$ is an $L$-dimensional subspace of $K^N$ and $C$ is any $(N-L) \times N$ matrix over $O_K$ of rank $1 \leq L < N$, viewed as a linear map $O_K^N \to O_K^{N-L}$, such that $V = \{ \bx \in K^N : C\bx = \bo \}$, let us write
$$H_{\inf}(C) = \prod_{v | \infty} \H_v(C)^{d_v/d},\ H_{\fin}(C) = \prod_{v \nmid \infty} H_v(C)^{d_v/d},$$
and then by Lemma 2.1 and Proposition 2.4 of \cite{liebendorf:1}, we have
\begin{equation}
\label{hqf2}
H_{\fin}(C) = \left[ O_K^{N-L} : C(O_K ^N) \right]^{-1/d}.
\end{equation}
This means that the definitions over $K$ and over $D$ are really analogous.

Now let $F(\bX,\bY) \in D[\bX,\bY]$ be a hermitian form in $2N$ variables with coefficients in $D$, so that $F(a \bx,\bwy) = \bar{a} F(\bx,\bwy)$ and $F(\bwy,\bx) = \overline{F(\bx,\bwy)}$ for each $a \in D$ and $\bx,\bwy \in D^N$. We also write $F(\bX)$ for $F(\bX,\bX)$, then $F(\bx) \in K$ for any $\bx \in D^N$. Let us also write $\Ff = (f_{ml})$ for the $N \times N$ coefficient matrix of $F$, then $f_{ml} = \overline{f_{lm}}$ for each $1 \leq l,m \leq N$, and $F(\bX,\bY) = \bX^t \Ff \bY$. In the same way as for quadratic and bilinear forms over $K$, we will talk about the height of the hermitian form $F$ over $D$, where by $H^{\O}(F)$ (respectively,  $H_{\inf}(F)$, $\Hf(F)$) we will always mean $H^{\O}(\Ff)$ (respectively,  $H_{\inf}(\Ff)$, $\Hf(\Ff)$), viewing $\Ff$ as a vector in $D^{N^2}$. We define the corresponding bilinear form $B$ over $K$ by taking the trace of $F$, i.e. $B([\bX],[\bY]) = \Tr(F(\bX,\bY))$. The associated quadratic form
\begin{equation}
\label{trace_form}
Q([\bX]) := B([\bX],[\bX])
\end{equation}
in $4N$ variables over $K$ is equal to $2F(\bX)$. Therefore $F(\bx) = 0$ for some $\bx \in D^N$ if and only if $Q([\bx]) = 0$. Write $\Bb$ for the $4N \times 4N$ symmetric matrix of $B$ over $K$, then each entry of $\Ff$ corresponds to a $4 \times 4$ block in $\Bb$. Specifically, if $f_{ml} = f_{ml}(0)+ f_{ml}(1)i+ f_{ml}(2)j+ f_{ml}(3)k \in D$, then the corresponding block in $\Bb$ is of the form
\begin{equation}
\label{trace_matrix}
\Bb(f_{ml}) := \left( \begin{matrix}
2f_{ml}(0)&2\alpha f_{ml}(1)&2\beta f_{ml}(2)&-2\alpha \beta f_{ml}(3)\\
-2\alpha f_{ml}(1)&-2\alpha f_{ml}(0)&-2\alpha \beta f_{ml}(3)&2\alpha \beta f_{ml}(2)\\
-2\beta f_{ml}(2)&2\alpha \beta f_{ml}(3)&-2\beta f_{ml}(0)&-2\alpha \beta f_{ml}(1)\\
2\alpha \beta f_{ml}(3)&-2\alpha \beta f_{ml}(2)&2\alpha \beta f_{ml}(1)&2\alpha \beta f_{ml}(0)
\end{matrix} \right),
\end{equation}
so $\Bb = (\Bb(f_{ml}))_{1 \leq m,l \leq N}$, and $Q(\bz) = \bz^t \Bb \bz$ for each $\bz \in K^{4N}$. As defined before, we will write $H(Q)$ (respectively,  $H_{\inf}(Q)$, $H_{\fin}(Q)$) for $H(\Bb)$ (respectively,  $H_{\inf}(\Bb)$, $H_{\fin}(\Bb)$), viewed as a vector in $K^{16N^2}$.
\smallskip

Finally, we define the constants that appear in our inequalities over quaternion algebras. Define a special order $\O_D$ in $D$:
\begin{equation}
\label{O_D}
\O_D = O_K + O_K i + O_K j + O_K k.
\end{equation}
For our fixed order $\O$, define
\begin{equation}
\label{cko}
c_{\O}(Z) = \min \left\{ h(a) : a \in K \text{ such that } a Z \cap \O^N \subset \O_D^N \right\},
\end{equation}
as well as
\begin{equation}
\label{zko}
z_{\O}(Z) = \min \left\{ h(a) h(a^{-1}) : a \in K \text{ such that } a Z \cap \O^N \subset \O_D^N \right\}.
\end{equation}
Let $\Delta_{\O}$ be the discriminant of the order $\O$, which is the ideal in $O_K$ generated by all the elements of the form
$$\det \left( \Tr (\omega_h \omega_n) \right)_{0 \leq h,n \leq 3} \in O_K,$$
where $\omega_0,\dots,\omega_3$ are in $\O$. Now the constants used in the statement of Theorem~\ref{main-1} are given by
\begin{equation}
\label{E3}
\E_3(D,\O,Z,d,L) = 2^{\frac{4L(d-2)+3}{2}} Ld\ s(\alpha,\beta) z_{\O}(Z) c_{\O}(Z)^{4Ld-1} \Nn(\Delta_{\O})^{\frac{L}{2}},
\end{equation}
where $\Nn$ stands for the norm from $K$ to $\que$, and
\begin{equation}
\label{E4}
\E_4(D,\O,Z,d,L) = \frac{c_{\O}(Z)}{2 \sqrt{2} Ld\ s(\alpha,\beta) z_{\O}(Z)}.
\end{equation}
We also define the constant that appears in the upper bound of Theorem \ref{main2}. Let
\begin{equation}
\label{O_const}
\MM(\O) := \max \left\{ \frac{\Nn(\Delta_{\O})^{1/2}}{\Nn(4 \alpha \beta)}, \frac{\Nn(4 \alpha \beta)}{\Nn(\Delta_{\O})^{1/2}} \right\},
\end{equation}
and define
\begin{equation}
\label{AK}
\A_{K,\O} (L,M,J,\alpha,\beta) = \frac{2^{\frac{9L+13}{2}} s(\alpha,\beta)^{9L+12}}{t(\alpha,\beta)^{\frac{9L+11}{2}}} \MM(\O)^{4(N-L)(9L+12)} T_K(L,M+2J+1),
\end{equation}
where the field constant $T_K(\ell,j)$ is defined in \eqref{TK} and $s(\alpha,\beta)$, $t(\alpha,\beta)$ are defined in (\ref{s_ab}). We are now ready to proceed.
\bigskip

\section{Counting points of bounded height}
\label{count}

Here we discuss counting estimates for the cardinality of sets of points of bounded height over number fields and quaternion algebras, as discussed above. In particular, we prove Theorems~\ref{cnt_module}, \ref{main-1} and~\ref{main-2}.

Our main tool is a basic counting mechanism for lattice points in cubes, which is a consequence of results of~\cite{me:classical} and~\cite{me:number}. Let us write $C_n(R)$ for the closed cube of side-length $2R$ centered at the origin in $\real^n$, i.e.
\begin{equation}
\label{cube}
C_n(R) = \left\{ \bx \in \real^n : \max_{1 \leq m \leq n} |x_m| \leq R \right\}.
\end{equation}

\begin{lem} \label{cnt_lem} Let $\Lambda \subset \real^N$ be a lattice of rank $L \leq N$ so that for every $\bo \neq \bx \in \Lambda$,
\begin{equation}
\label{max_c}
|\bx| := \max_{1 \leq n \leq N} |x_n| \geq c
\end{equation}
for some $c \in \real_{>0}$ independent of $\bx$. Then for any $R \in \real_{>0}$,
\begin{equation}
\label{lattice_up}
\left| \Lambda \cap C_N(R) \right| \leq \left\{ \begin{array}{ll}
\left( \frac{2Rc^{N-1}}{\det(\Lambda)} + 1 \right) \left( \frac{2R}{c} + 1 \right)^{N-1} & \mbox{if $L=N$} \\
\left( \frac{2R}{c} + 1 \right)^{N-1} & \mbox{if $L<N$} \\
\left( \frac{2 \binom{N}{L}^{1/2} R}{\det(\Lambda)} + 1 \right) \left( 2R + 1 \right)^{L-1} & \mbox{if $\Lambda \subseteq \zed^N$}
\end{array}
\right.
\end{equation}
In addition, if $R \geq \frac{L}{2} \max \left\{ \frac{\det(\Lambda)}{c^{L-1}}, c \right\}$, then
\begin{equation}
\label{lattice_low}
\left| \Lambda \cap C_N(R) \right| \geq \left( \frac{2Rc^{L-1}}{L\det(\Lambda)} - 1 \right) \left( \frac{2R}{Lc} - 1 \right)^{L-1}.
\end{equation}
\end{lem}

\proof
We start by obtaining the upper bound of~\eqref{lattice_up}. If $L=N$, then~\eqref{lattice_up} follows from Lemma~2.1 of~\cite{me:number}. Assume that $L < N$ and let
$$V = \left\{ \bx \in \real^N : \bx \cdot \bwy = 0\ \forall\ \bwy \in \Lambda \right\}$$
be the $(N-L)$-dimensional subspace of $\real^N$ orthogonal to $\Lambda$. Let $\Lambda' \subseteq V$ be a full-rank lattice in $V$ spanned by an orthogonal basis of unit vectors, then $\det(\Lambda') = 1$ and a shortest nonzero vector in $\Lambda'$ has norm~1. Now let $R \in \real_{>0}$ and let $R_* \geq \max \left\{ 2\sqrt{N} R, \frac{c\sqrt{N}}{\min \{ |\bx| : \bo \neq \bx \in \Lambda' \}} \right\}$, and define
$$\Lambda''(R_*) = \Lambda \oplus R_* \Lambda'.$$
Notice that every $\bx \in \Lambda''(R_*)$ is of the form $\bx = \bx_1 + R_*\bx_2$ for some $\bx_1 \in \Lambda$, $\bx_2 \in \Lambda'$ and $\|\bx\|^2 = \|\bx_1\|^2 + R_*^2\|\bx_2\|^2$ since $\bx_1$ and $\bx_2$ are orthogonal. Therefore
\begin{equation}
\label{R_bnd}
|\bx| \geq \frac{1}{\sqrt{N}} \|\bx\| = \frac{1}{\sqrt{N}} \sqrt{\|\bx_1\|^2 + R_*^2 \|\bx_2\|^2} \geq \max\{ 2R, c \},
\end{equation}
which in particular means that if $\bx \in \Lambda''(R_*) \cap C_N(R)$, then $\bx_2 = \bo$ and so $\bx \in \Lambda$. Hence $|\Lambda \cap C_N(R)| = |\Lambda''(R_*) \cap C_N(R)|$ and
\begin{equation}
\label{det_LO-1}
\det(\Lambda''(R_*)) = R_*^{N-L} \det(\Lambda) \geq (2NR)^{N-L} \det(\Lambda).
\end{equation}
Since rank of $\Lambda''(R_*)$ is $N$, combining~\eqref{R_bnd} and~\eqref{det_LO-1} with Lemma~2.1 of~\cite{me:number} produces the bound
$$|\Lambda \cap C_N(R)| \leq \left( \frac{2Rc^{N-1}}{R_*^{N-L} \det(\Lambda)} + 1 \right) \left( \frac{2R}{c} + 1 \right)^{N-1},$$
and the bound of~\eqref{lattice_up} in case $L<N$ follows by taking the limit as $R_* \to \infty$.

Next, following~\cite{me:number}, let $X$ be a basis matrix for $\Lambda$ and for each $I \subset \{ 1,\dots,N \}$ with $|I| = L$ write $X_I$ for the $L \times L$ submatrix of $X$ whose columns are indexed by the elements of $I$. Let $J \subset \{ 1,\dots,N \}$ with $|J| = L$ be such that
$$|\det(X_J)| = \max_{|I|=L} |\det(X_I)|,$$
and let $\Omega$ be the lattice of full rank in $\real^{L}$ spanned over $\zed$ by the column vectors of~$X_J$. Then $\det(\Omega) = |\det(X_J)|$ is maximum of absolute values of Grassmann coordinates of $\Lambda$, and Cauchy-Binet formula (see, for instance (18) of~\cite{me:number}) implies that
\begin{equation}
\label{det_LO}
\det(\Omega) \leq \det(\Lambda) \leq \binom{N}{L}^{1/2} \det(\Omega).
\end{equation}
The bound of~\eqref{lattice_up} in case $\Lambda \subseteq \zed^N$ follows by combining~\eqref{det_LO} with Theorem~4.2 of~\cite{me:classical}.

Now we derive the lower bound of~\eqref{lattice_low}. By Corollary 1 on p. 13 of \cite{cassels_geom}, it is possible to select a basis for $\Omega$ such that the basis matrix $A$ is upper triangular, all of its nonzero entries are positive, and the maximum entry of each row occurs on the diagonal. By~\eqref{max_c} above, each of these entries is at least $c$, since each column of $A$ is a linear combination of columns of $X_J$. The lattice $\Omega$ now satisfies the conditions of Lemma~2.1 of~\cite{me:number}, and so if $2R \geq \max \left\{ \frac{\det(\Omega)}{c^{L-1}}, c \right\}$, then
\begin{equation}
\label{O_lattice_low}
\left| \Omega \cap C_L(R) \right| \geq \left( \frac{2Rc^{L-1}}{\det(\Omega)} - 1 \right) \left( \frac{2R}{c} - 1 \right)^{L-1},
\end{equation}
where the condition on $R$ simply ensures that every term in the product on the right hand side of the inequality is positive. Now Theorem 4.3 (equation~(31)) of~\cite{me:classical} implies that
\begin{equation}
\label{O_bounds-1}
\left| \Lambda \cap C_{N}(R) \right| \geq \left| \Omega \cap C_{L}\left( \frac{R}{L} \right) \right|,
\end{equation}
and combining this observation with~\eqref{det_LO} and~\eqref{O_lattice_low}, we obtain~\eqref{lattice_low}.
\endproof

We now use Lemma~\ref{cnt_lem} to prove Theorem~\ref{cnt_module}, producing an estimate on the number of points of bounded height in a fixed torsion-free $O_K$-module for an arbitrary number field~$K$.

\proof[Proof of Theorem~\ref{cnt_module}]  Let all the notation be as in the statement of the theorem. Since $\M \subset K^N$, it must be torsion-free, hence projective. By the structure theorem for finitely generated projective modules over Dedekind domains (see, for instance~\cite{lang}),
$$\M = \left\{ \sum_{n=1}^L \beta_n \bwy_n : \bwy_n \in O_K^N,\ \beta_n \in \I_n \right\}$$
for some $O_K$-fractional ideals $\I_1,\dots,\I_L$ in $K$. By Proposition~13 on p.66 of~\cite{lang}, the discriminant of $\M$ is then
\begin{equation}
\label{module_disc}
\D_K(\M) := \D_K^L \prod_{n=1}^L \Nn(\I_n)^2,
\end{equation}
where $\Nn(\I_n)$ is the norm of the fractional ideal $\I_n$. Define $\UU_K(\M)$, a fractional $O_K$-ideal in $K$, to be
\begin{equation}
\label{UKM}
\UU_K(\M) = \left\{ \alpha \in K : \alpha\M \subseteq O_K^N \right\},
\end{equation}
then
$$c_K(\M) = \min \{ h(\alpha) : \alpha \in \UU_K(\M) \}.$$
Let 
$$\sigma_1,\dots,\sigma_{r_1},\tau_1,\dots,\tau_{r_2},\dots,\tau_{2r_2}$$
be the embeddings of $K$ into $\cee$ with $\sigma_1,\dots,\sigma_{r_1}$ being the real embeddings and $\tau_n,\tau_{r_2+n} = \bar{\tau}_n$ for each $1 \leq n \leq r_2$ being the pairs of complex conjugate embeddings. For each $\alpha \in K$ and each complex embedding $\tau_n$, write $\tau_{n1}(\alpha) = \Re(\tau_n(\alpha))$ and $\tau_{n2}(\alpha) = \Im(\tau_n(\alpha))$, where $\Re$ and $\Im$ stand respectively for real and imaginary parts of a complex number. Then $d=r_1+2r_2$, and we define an embedding
$$\sigma^N = (\sigma_1^N,\dots,\sigma_{r_1}^N, \tau_{11}^N,\tau_{12}^N,\dots,\tau_{r_21}^N,\tau_{r_22}^N) : K^N \to \real^{Nd}.$$
Let $\alpha \in \UU_K(\M)$. Since $\alpha \bx \in O_K^N$ for every $\bx \in \M$, we have
\begin{eqnarray*}
\max \{ |\sigma_1(\alpha x_n)|,\dots,|\sigma_{r_1}(\alpha x_n)|,|\tau_{11}(\alpha x_n)|,|\tau_{12}(\alpha x_n)|,\dots,|\tau_{r_2 1}(\alpha x_n)|,|\tau_{r_2 2}(\alpha x_n)| \} \\
\geq  \frac{1}{\sqrt{2}},\ \ \ \ \ \ \ \ \ \ \ \ \ \ \ \ \ \ \ \ \ \ \ \ \ \ \ \ \ \ \ \ \ \ \ \ \ \ \ \ \ \ \ \ \ \ \ \ \ \ \ \ \ \ \ \ \ \ \ \ \ \ \ \ \ \ \ \ \ \ \ \ \ \ \ \ \ \ \ \ \ \ \ \ \ \ \ \ \ \ \ \ \ \ \ \ \ \ 
\end{eqnarray*}
for every $1 \leq n \leq N$, as indicated in~\cite{me:number}, and therefore
\begin{eqnarray*}
& & \max \{ |\sigma_1(x_n)|,\dots,|\sigma_{r_1}(x_n)|,|\tau_{11}(x_n)|,|\tau_{12}(x_n)|,\dots,|\tau_{r_2 1}(x_n)|,|\tau_{r_2 2}(x_n)| \} \\
& \geq & \frac{1}{\sqrt{2}} \max \{ |\sigma_1(\alpha)|,\dots,|\sigma_{r_1}(\alpha)|,|\tau_{11}(\alpha)|,|\tau_{12}(\alpha)|,\dots,|\tau_{r_2 1}(\alpha)|,|\tau_{r_2 2}(\alpha)| \}^{-1} \\
& \geq & \frac{1}{\sqrt{2}} \prod_{l=1}^{r_1} \max \{ 1, |\sigma_l(\alpha)| \}^{-1} \times \prod_{m=1}^{r_2} \max \{ 1, |\tau_m(\alpha)| \}^{-1} \\
& \geq & \frac{1}{\sqrt{2}} h(\alpha)^{-1}.
\end{eqnarray*}
Since the choice of $\alpha \in \UU_K(\M)$ was arbitrary, we can pick such an $\alpha$ with $h(\alpha) = c_K(\M)$, and so
\begin{eqnarray}
\label{M_c}
& & \max \{ |\sigma_1(x_n)|,\dots,|\sigma_{r_1}(x_n)|,|\tau_{11}(x_n)|,|\tau_{12}(x_n)|,\dots,|\tau_{r_2 1}(x_n)|,|\tau_{r_2 2}(x_n)| \} \nonumber \\
& & \geq \frac{1}{\sqrt{2}} c_K(\M)^{-1}
\end{eqnarray}
for every $1 \leq n \leq N$, $\bx \in \M$. Notice that $\Lambda_K(\M) := \sigma^N(\M)$ is a lattice of rank $Ld$ in $\real^{Nd}$, and a direct adaptation of Lemma~2 on p.115 of~\cite{lang} implies that the determinant of $\Lambda_K(\M)$ is
\begin{equation}
\label{det_lkm}
\det(\Lambda_K(\M)) = 2^{-Lr_2} |\D_K(\M)|^{\frac{L}{2}} = 2^{-Lr_2} |\D_K|^{\frac{L}{2}} \prod_{n=1}^L \Nn(\I_n),
\end{equation}
where the last identity follows by~\eqref{module_disc} above. Combining~\eqref{M_c} and~\eqref{det_lkm} with Lemma~\ref{cnt_lem}, we see that the cardinality of the set $\Lambda_K(\M) \cap C_{Nd}(R)$ is
\begin{equation}
\label{mn_lat_bnd}
\geq \left( \frac{R}{2^{\frac{Lr_1-3}{2}} Ld\ c_K(\M)^{Ld-1} |\D_K(\M)|^{\frac{L}{2}}} - 1 \right) \left( \frac{2^{\frac{3}{2}} c_K(\M)R}{Ld} - 1 \right)^{Ld-1}.
\end{equation}
For any $\alpha \in \UU_K(\M)$, $\alpha \bx \in O_K^N$ for every $\bx \in \M$, and so
$$\prod_{v \in M(K), v \nmid \infty} \max \{1, |\alpha x_1|_v, \dots, |\alpha x_N|_v \} = 1,$$
and so
\begin{eqnarray*}
h(\alpha \bx)^d & = & \prod_{m=1}^{r_1} \max \{ 1, |\sigma_m(\alpha x_1)|,\dots,|\sigma_m(\alpha x_{N})| \} \times \\
& & \times \prod_{n=1}^{r_2} \max \{ 1, \tau_{n1}(\alpha x_1)^2+\tau_{n2}(\alpha x_1)^2,\dots,\tau_{n1}(\alpha x_N)^2+\tau_{n2}(\alpha x_N)^2 \} \\
& \leq & h(\alpha)^d |\sigma^N(\bx)|^d,
\end{eqnarray*}
where $|\bwy| = \max_{1 \leq n \leq Nd} |y_n|$ for each vector $\bwy \in \real^{Nd}$, and so
$$h(\bx) = h(\alpha^{-1} (\alpha \bx)) \leq h(\alpha^{-1}) h(\alpha \bx) \leq h(\alpha^{-1}) h(\alpha) |\sigma^N(\bx)|.$$
Notice that
$$z_K(\M) = \min \left\{ h(\alpha) h(\alpha^{-1}) : \alpha \in \UU_K(\M) \right\},$$
and choose $\alpha$ with $h(\alpha) h(\alpha^{-1}) = z_K(\M)$, then
\begin{equation}
\label{bnd_x_z}
h(\bx) \leq z_K(\M) |\sigma^N(\bx)|
\end{equation}
for every $\bx \in \M$. Therefore 
$$\Lambda_K(\M) \cap C_{Nd}(R) \subseteq \sigma^N(S_{\M}(z_K(\M) R)).$$
Combining this observation with~\eqref{mn_lat_bnd} yields~\eqref{mn_low_bnd}.
\endproof

We will now apply the bound of Theorem~\ref{cnt_module} to obtain a lower bound on the number of points of bounded height in a right $D$-vector space which are integral over a fixed order $\O$ in $D$.

\proof[Proof of Theorem~\ref{main-1}] Define $Z_{\O} = Z \cap \O^N$ and let $\M_Z = [Z_{\O}] \subset K^{4N}$, which is an $O_K$-module such that $\M_Z \otimes_K K \cong K^{4L}$. Suppose that $\bwy \in \M_Z$ satisfies $h(\bwy) \leq R$, then $\bx := [\bwy]^{-1} \in Z_{\O}$ and
$$h(\bx) \leq 2s(\alpha,\beta) h(\bwy) \leq 2s(\alpha,\beta) R,$$
by Lemma~3.1 of~\cite{quaternion}. Therefore
\begin{equation}
\label{bnd_M1}
\left| S_{D,N}(Z,\O,R) \right| \geq \left| \left\{ \bwy \in \M_Z : h(\bwy) \leq \frac{R}{2s(\alpha,\beta)} \right\} \right|.
\end{equation}
Therefore we can apply Theorem~\ref{cnt_module} to $\M_Z$, obtaining a lower bound on the number of points of bounded height in~$\M_Z$. To derive~\eqref{count_ZO} from this bound, we need to relate invariants of $\M_Z$ which appear in~\eqref{mn_low_bnd} to corresponding invariants of $Z_{\O}$ and then apply the height comparison lemmas of~\cite{quaternion}.

Let $\Lambda_K(\M_Z) = \sigma^N(\M_Z)$, as in the proof of Theorem~\ref{cnt_module} above. Then Lemma~3.2 of~\cite{liebendorf:1} (also see the proof of Lemma~3.5 of~\cite{quaternion}) combined with equation \eqref{det_lkm} above asserts that
\begin{equation}
\label{det_M_Z}
|\D_K(\M_Z)|^{\frac{4L}{2}} = \det(\Lambda_K(\M_Z)) = \left( \sqrt{\Nn(\Delta_{\O})}/16 \right)^L H^{\O}(Z)^{4d}.
\end{equation}
Also notice that $aZ_{\O} \subseteq \O_D^N$ for some $a \in K$ if and only $a\M_Z \subseteq O_K^{4N}$, which means that $c_{\O}(Z) = c_K(\M_Z)$ and $z_{\O}(Z) = z_K(\M_Z)$, where $c_K(\M_Z)$ and $z_K(\M_Z)$ are defined as in~\eqref{ckm} and~\eqref{zkm} above. Now combining~\eqref{mn_low_bnd} with~\eqref{bnd_M1} and~\eqref{det_M_Z}, we see that $\left| S_{D,N}(Z,\O,R) \right| \geq$
\begin{eqnarray}
\label{asmpt_bnd}
& \geq & \left( \frac{R}{\E_3(D,\O,Z,d,L) H^{\O}(Z)^{4d}} - 1 \right) \left( \E_4(D,\O,Z,d,L) R - 1 \right)^{4Ld-1} \nonumber \\
& = & \E'_3(D,\O,Z,d,L) \frac{R^{4Ld}}{H^{\O}(Z)^{4d}} + O(R^{4Ld-1}),
\end{eqnarray}
where
$$\E'_3(D,\O,Z,d,L) = \left( 2^{4L(2d-1)} (Ld\ s(\alpha,\beta) z_{\O}(Z))^{4Ld} \Nn(\Delta_{\O})^{\frac{L}{2}} \right)^{-1}.$$
This finishes the proof.
\endproof

Finally, we apply the counting estimate of~\cite{loher_masser} over number fields to prove Theorem~\ref{main-2}.

\proof[Proof of Theorem~\ref{main-2}] Since $[\ ] : D \to K^4$ is a vector space isomorphism,
$$\left| S_{D,N}(R) \right| = \left| [S_{D,N}(R)] \right|.$$
Now Lemma~3.1 (or, more precisely, inequality (18)) of \cite{quaternion} guarantee that for every $\bx \in D^N$,
\begin{equation}
\label{ht_3.1}
t(\alpha,\beta) h([\bx]) \leq h(\bx) \leq s(\alpha,\beta) h([\bx]),
\end{equation}
and hence
$$[S_{D,N}(R)] \subseteq S_{K,4N}(R/t(\alpha,\beta)) := \left\{ \bwy \in K^{4N} : h(\bwy) \leq \frac{R}{t(\alpha,\beta)} \right\}.$$
An upper bound on cardinality of the set $S_{K,4N}(R/t(\alpha,\beta))$ follows from Theorem~4 of \cite{loher_masser}:
\begin{equation}
\label{count_up}
\left| S_{K,4N} \left( R/t(\alpha,\beta) \right) \right| \leq (1088 d \log d)^{4N} \left( \frac{R}{t(\alpha,\beta)} \right)^{(4N+1)d}.
\end{equation}
\endproof

\begin{rem} \label{cnt}
On the other hand, \eqref{ht_3.1} implies that
$$S_{K,4N}(R/s(\alpha,\beta)) := \left\{ \bwy \in K^{4N} : h(\bwy) \leq \frac{R}{s(\alpha,\beta)} \right\} \subseteq [S_{D,N}(R)].$$
Equation (1.5) of \cite{schmidt_northcott-1} implies that
\begin{equation}
\label{count_low}
\left| S_{K,4N} \left( R/s(\alpha,\beta) \right) \right| \gg_{K,N} \left( \frac{R}{s(\alpha,\beta)} \right)^{4N+1}.
\end{equation}
Then \eqref{count_bnd} follows by combining \eqref{count_up} with \eqref{count_low}. In fact, as long as we have {\it any} upper or lower bounds on the number of points of bounded height over $K$, we can ``transfer" them to obtain analogous bounds for the number of points of bounded height over $D$.
\end{rem}
\bigskip

\appendix
\section{Further counting estimates over global fields}
\label{S-unit}

Here we show some further applications of Lemma~\ref{cnt_lem}, obtaining estimates on the number of $S$-units of bounded height in an arbitrary number field as well as number of rational functions of bounded height supported on a given curve over a fixed finite field.

We start with the number field situation. Let $K$ be any number field, and write $S_{\infty}$ for the set of all archimedean places of $K$. Let $S_1$ be a finite (possibly empty) set of non-archimedean places of $K$, and let $S = S_{\infty} \cup S_1$. The group of $S$-units of $K$ is
$$O^*_S = \left\{ a \in K : |a|_v = 1\ \forall\ v \notin S \right\}.$$
Define the logarithmic $S$-height function on $K^{\times}$ by
\begin{equation}
\label{s-height}
H_S(a) = \max_{v \in S} \{ \left| \log |a|_v \right|, \left| \log |a^{-1}|_v \right| \},
\end{equation}
i.e., $H_S(a)$ measures the extent of divisibility of numerator and denominator of $a$ at the places in $S$, and let 
\begin{equation}
\label{HSK}
H_{S,K} = \min \{ H_S(a) : a \in O_S^* \setminus \mu_K \} > 0,
\end{equation}
where $\mu_K$ is the group of roots of unity in $K$. Let $d=[K:\que]$, $h_K$ be the class number and $R_K$ the regulator of $K$. 

We employ here the standard logarithmic lattice construction used in the proof of Dirichlet's Unit Theorem (see, for instance, p.104 of~\cite{lang} and pp.575--578 of~\cite{tsfasman}). Let $n = |S| = d + t$, where $t=|S_1|$, and define the map $\varphi_S : O_S^* \to \real^n$ by
$$\varphi_S(a) = (\log |a|_v)_{v \in S}.$$
Then $\Ker \varphi_S = \mu_K$ and $L_S := \varphi(O_S^*)$ is a lattice of rank $(n-1)$ in $\real^n$, which is contained in the hyperplane $V = \{ \bx \in \real^n : \sum_{m=1}^n  x_m = 0 \}$, and so $L_S$ is a lattice of full rank in~$V$. The $S$-regulator of $K$ is defined to be
$$R_{S,K} := \det L_S,$$
which is just $R_K$ if $S_1 = \emptyset$. If $S_1 \neq \emptyset$, let $\pp_1,\dots,\pp_t$ be the prime ideals in $K$ corresponding to the places in~$S_1$, and let $P$ be the largest rational prime lying below these prime ideals. In Lemma~3 of~\cite{bugeaud}, the following bounds on $R_{S,K}$ are produced (see also Lemma~3 of~\cite{bg} and Proposition~5.4.7 of~\cite{tsfasman}):
\begin{equation}
\label{RS_up}
R_{S,K} \leq R_K h_K \prod_{m=1}^t \log \Nn(\pp_m) \leq R_K h_K (d \log^* P)^t
\end{equation}
and
\begin{equation}
\label{RS_low}
R_{S,K} \geq R_K \prod_{m=1}^t \log \Nn(\pp_m) \geq 0.2052 (\log 2)^d (\log^* P),
\end{equation}
where $\log^* P = \max \{ \log P, 1\}$. Observe also that for any $\bx \in L_S \setminus \{\bo\}$,
\begin{equation}
\label{mx_const}
|\bx| = \max_{1 \leq m \leq n} |x_m| \geq H_{S,K} > 0.
\end{equation}
We are now ready to state and prove our estimate.

\begin{lem} \label{sunit} Let $B \in \real_{>0}$ and let
$$O^*_S(B) = \left\{ a \in O^*_S : H_S(a) \leq B \right\}.$$
Then, with notation as above,
\begin{eqnarray}
\label{sunit_bound}
\omega_K \left( \frac{2BH_{S,K}^{n-2}}{(n-1)R_{S,K}} - 1 \right) \left( \frac{2B}{(n-1)H_{S,K}} - 1 \right)^{n-2} \nonumber \\ 
 \leq |O^*_S(B)| \leq \omega_K \left( \frac{2B}{H_{S,K}} + 1 \right)^{n-1},
\end{eqnarray}
where $\omega_K = |\mu_K|$; the lower bound of~\eqref{sunit_bound} holds for $B \geq \frac{n-1}{2} \max \left\{ \frac{R_{S,K}}{H_{S,K}^{n-2}}, H_{S,K} \right\}$.
\end{lem}

\proof Given a positive real number $B$, let $C_n(B)$ be as in~\eqref{cube}. It is then an easy observation that $O_S^*(B) = \varphi_S^{-1}(C_n(B) \cap L_S)$. Notice that for each $\bx \in L_S$, $|\varphi_S^{-1}(\bx)| = \omega_K$, therefore
\begin{equation}
\label{in_out}
|O^*_S(B)| = \omega_K\ |C_n(B) \cap L_S|,
\end{equation}
and~\eqref{sunit_bound} follows by combining~\eqref{mx_const} and~\eqref{in_out} with Lemma~\ref{cnt_lem}.
\endproof

\begin{rem} \label{sunit_rem} Inequalities~\eqref{RS_up} and~\eqref{RS_low} can now be used to make estimates of Lemma~\ref{sunit} more explicit, if necessary. Comparable asymptotic estimates on the number of units and $S$-units of bounded height (with somewhat different heights used) were previously obtained in~\cite{everest} (see also Theorem~5.2 on p.70 of~\cite{lang_dioph}) and~\cite{fuchs} (Lemma 1). In contrast, our estimates are explicit upper and lower bounds.
\end{rem}
\smallskip

Next we discuss an analogous construction over function fields, following pp.578--581 of~\cite{tsfasman}. Let $q$ be a prime power and let $\Ff_q$ be the finite field with $q$ elements. Let $X$ be a smooth projective curve defined over $\Ff_q$, and let $K=\Ff_q(X)$ be the field of rational functions on $X$ over $\Ff_q$. For every $f \in K^{\times}$, we write $\Supp(f)$ for the support of $f$, i.e., the set of all points at which $X$ has zeros or poles. Let $X(\Ff_q)$ be the set of points of $X$ which are rational over $\Ff_q$. Let $\Pp \subseteq X(\Ff_q)$, and define
$$O_{\Pp}^* = \{ f \in K^{\times} : \Supp(f) \subseteq \Pp \}$$
to be the group of all rational functions in $K$ supported on $\Pp$. Let $n = |\Pp|$, say $\Pp = \{ p_1,\dots,p_n\}$, and write $a_m(f) \in \zed$ for the order of zero or pole that $f \in K^{\times}$ has at $p_m \in \Pp$. We can define the $\Pp$-height on $K^{\times}$ by
\begin{equation}
\label{p-height}
H_{\Pp}(f) = \max_{1 \leq m \leq n} |a_m(f)|,
\end{equation}
which is a direct function-field analogue of the $S$-height function defined in~\eqref{s-height} above. The principal divisor of any $f \in O_{\Pp}^*$ is
$$\div(f) = a_1(f)p_1+\dots+a_n(f)p_n,$$
so that $\sum_{m=1}^n a_m(f) = 0$. We can then define a map $\varphi_{\Pp} : O_{\Pp}^* \to \real^n$ by
$$\varphi_{\Pp}(f) = (a_1(f),\dots,a_n(f)),$$
and so $\Ker(\varphi_{\Pp}) = \Ff_q^{\times}$ and $L_{\Pp} :=\varphi_{\Pp}(O_{\Pp}^*)$ is a finite-index sublattice of the root lattice
$$A_{n-1} = \left\{ \bx \in \zed^n : \sum_{m=1}^n x_m = 0 \right\},$$
which has determinant $=\sqrt{n}$. We need some more notation to give a formula for the determinant of $L_{\Pp}$, following~\cite{tsfasman}. Let $\Div^0(X)$ be the group of divisors of degree 0 on $X$ and $P(X)$ the subgroup of principal divisors, then $J(X)=\Div^0(X)/P(X)$ is the Jacobian of $X$, and we write $J_X(\Ff_q)$ for the set of $\Ff_q$-rational points on the Jacobian. Let also $\Div^0_{\Pp}(X) \subset \Div^0(X)$ be the subgroup of degree 0 divisors supported on $\Pp$ and $P_{\Pp}(X) = \Div^0_{\Pp}(X) \cap P(X)$. Define the restricted $\Pp$-Jacobian to be $J_{X,\Pp} := \Div^0_{\Pp}(X) / P_{\Pp}(X)$, then Theorem~5.4.9 of~\cite{tsfasman} states that
\begin{equation}
\label{det_P}
\det(L_{\Pp}) = \det(A_{n-1}) |A_{n-1} : L_{\Pp}| = \sqrt{n} \left| J_{X,\Pp} \right|,
\end{equation}
and so
\begin{equation}
\label{det_P_1}
\sqrt{n} \leq \det(L_{\Pp}) \leq \sqrt{n} |J_X(\Ff_q)| \leq \sqrt{n} \left( 1+q+\frac{|X(\Ff_q)|-q-1}{g} \right)^g,
\end{equation}
where $g$ is the genus of $X$. Further, the same theorem guarantees that for every $\bx \in L_{\Pp} \setminus \{\bo\}$,
\begin{equation}
\label{max-x-curve}
|\bx| \geq \max \left\{ 1, \frac{1}{n} \sqrt{\frac{2|X(\Ff_q)|}{q+1}} \right\}.
\end{equation}
We are now ready to state and prove the function-field analogue of Lemma~\ref{sunit}.

\begin{lem} \label{p-count} Let $B \in \real_{>0}$ and let
$$O^*_{\Pp}(B) = \left\{ f \in O^*_{\Pp} : H_{\Pp}(f) \leq B \right\}.$$
Then, with notation as above,
\begin{eqnarray}
\label{p_bound}
(q-1) \left( \frac{2B}{(n-1) \sqrt{n} \left| J_{X,\Pp} \right|} - 1 \right) \left( \frac{2B}{(n-1)} - 1 \right)^{n-2} \nonumber \\ 
 \leq |O^*_{\Pp}(B)| \leq (q-1) \left( \frac{2B}{\left| J_{X,\Pp} \right|} + 1 \right) \left( 2B + 1 \right)^{n-2}
\end{eqnarray}
where the lower bound of~\eqref{sunit_bound} holds for $B \geq \frac{(n-1) \sqrt{n} \left| J_{X,\Pp} \right|}{2}$.
\end{lem}

\proof Given a positive real number $B$, let $C_n(B)$ be as in~\eqref{cube}. It is then an easy observation that $O_{\Pp}^*(B) = \varphi_{\Pp}^{-1}(C_n(B) \cap L_{\Pp})$. Notice that for each $\bx \in L_{\Pp}$, $|\varphi_{\Pp}^{-1}(\bx)| = q-1$, therefore
\begin{equation}
\label{in_out-1}
|O^*_{\Pp}(B)| = (q-1)\ |C_n(B) \cap L_{\Pp}|,
\end{equation}
and~\eqref{p_bound} follows by combining~\eqref{max-x-curve} and~\eqref{in_out-1} with Lemma~\ref{cnt_lem}.
\endproof

\begin{rem} \label{pcount_rem} Formulas~\eqref{det_P} and~\eqref{det_P_1} can be used to make estimates of Lemma~\ref{p-count} more explicit, if necessary.
\end{rem}

\bigskip

\section{Points of small height}
\label{points}

Classical Diophantine results on existence of points of bounded height on linear and quadratic spaces, such as Siegel's lemma and Cassels' theorem, have enjoyed much attention, including a number of papers by various authors in the recent years. In particular, some of the recent work has been devoted to extending these results to the non-commutative situation (see \cite{liebendorf:1}, \cite{liebendorf:2}, \cite{liebendorf:3}, \cite{watanabe}, \cite{quaternion}, and others). On the other hand, the non-commutative situation presents various obstacles that do not exist over fields, which makes it difficult to push the theory much further even over quaternion algebras. It is however possible to ``transfer" some of the existent results in the context of number fields to quaternion algebras, using appropriate height comparison inequalities. Here we demonstrate this transfer principle on several examples in the hope that it can also prove to be useful in a variety of other situations. As above, let $K$ be a totally real number field of degree $d$ over $\que$, and let $D = \binom{\alpha,\beta}{K}$ be a positive definite quaternion algebra over~$K$. Suppose we want to prove the existence of a nonzero point  $\bx \in D^N$ of explicitly bounded height which would satisfy a certain set of algebraic conditions. We suggest the use of the following basic method:
\smallskip

{\it Suppose we know that there exists a point $\bwy \in K^{4N}$ of bounded height such that $[\bwy]^{-1} \in D^N$ satisfies the desired algebraic conditions. Use the height comparison lemmas developed in Section~3 of~\cite{quaternion} to produce the necessary bounds on the height of $\bx := [\bwy]^{-1} \in D^N$.}
\smallskip

In other words, the results on points of bounded height over $D$ can be obtained by ``transferring" the analogous results over $K$ with the use of height comparison inequalities. The first instance of this method at work has been demonstrated in \cite{quaternion}, where a result on existence of a small-height basis for a hermitian space over $D$ consisting of zeros of the corresponding quadratic form has been obtained by the transfer of an analogous result over $K$, due to Vaaler~\cite{vaaler:smallzeros2}. We also used this same method above to derive Theorem~\ref{main-1} from Theorem~\ref{cnt_module}. Here we take this principle further, proving the analogues of some recent results of \cite{quad_zero} and \cite{null}. 

\begin{thm} \label{main1} Let $D = \binom{\alpha,\beta}{K}$ be a positive definite quaternion algebra over a totally real number field $K$, where $\alpha, \beta$ are totally negative algebraic integers in $K$. Let $\O$ be an order in $D$. Let $N \geq 2$ be an integer, and let $Z \subseteq D^N$ be an $L$-dimensional right $D$-subspace, $1 \leq L \leq N$. Let $U_1,\dots,U_M \subset D^N$ be proper right $D$-subspaces, let 
$$G_1(\bX,\bY),\dots,G_J(\bX,\bY) \in D[\bX,\bY]$$
be a hermitian forms in $2N$ variables, and let
\begin{equation}
\label{W_l}
W_l = \{ \bx \in D^N : G_l(\bx) := G_l(\bx,\bx) = 0 \}
\end{equation}
for each $1 \leq l \leq J$. Suppose that $Z \not\subseteq \left( \bigcup_{m=1}^M U_m \right) \left( \bigcup_{l=1}^J W_l \right)$. Then there exists a basis 
\begin{equation}
\label{basis_1}
\bwy_1,\dots,\bwy_L \in Z \setminus \left( \left( \bigcup_{m=1}^M U_m \right) \left( \bigcup_{l=1}^J W_l \right) \right)
\end{equation}
for $Z$ over $D$, such that
\begin{eqnarray}
\label{mn1}
h(\bwy_1) & \leq & h(\bwy_2) \leq \dots \leq h(\bwy_L) \nonumber \\
& \leq & 4 L (M+2J+1)^{\frac{1}{d}} |\D_K|^{\frac{L+1}{2d}} s(\alpha,\beta) \MM(\O)^{4(N-L)} H^{\O}(Z)^4.
\end{eqnarray}
\end{thm}

\proof
For an $L$-dimensional right $D$-subspace $Z \subseteq D^N$, $[Z]$ is a $4L$-dimensional subspace of $K^{4N}$. Recall from the definitions in Section~\ref{heights} that for a hermitian form $F(\bX,\bY) \in D[\bX,\bY]$ in $2N$ variables, its associated trace form
\begin{equation}
\label{trace}
Q_F([\bX]) = \Tr(F(\bX)) = F(\bX) + \overline{F(\bX)},
\end{equation}
which is a quadratic form in $4N$ variables over $K$, and $F(\bx) = 0$ for some $\bx \in D^N$ if and only if $Q_F([\bx])=0$. Then for each $W_l$ as in \eqref{W_l}, define
$$[W_l] := \{ [\bwy] \in K^{4N} : \bwy \in D^N,\ G_l(\bwy) = 0 \} = \{ \bx \in K^{4N} : Q_{G_l}(\bx) = 0 \}.$$

Now Theorem~A.1 of \cite{quad_zero} guarantees that there exists a basis $\bx_1,\dots,\bx_{4L}$ for $[Z]$ over $K$ such that 
$$\bx_1,\dots,\bx_{4L} \in [Z] \setminus \left( \left( \bigcup_{m=1}^M [U_m] \right) \left( \bigcup_{l=1}^J [W_l] \right) \right)$$
and
\begin{equation}
\label{h_order}
H(\bx_1) \leq H(\bx_2) \leq \dots \leq H(\bx_{4L}),\ h(\bx_1) \leq h(\bx_2) \leq \dots \leq h(\bx_{4L}),
\end{equation}
and for each $1 \leq n \leq 4L$,
\begin{equation}
\label{gen_bnd_1}
H(\bx_n) \leq h(\bx_n) \leq  2 L (M+2J+1)^{\frac{1}{d}} |\D_K|^{\frac{L+1}{2d}} H([Z]).
\end{equation}
Moreover, these vectors can be taken with coordinates in $O_K$. Notice that there exist
$$1=l_1 < l_2 < \dots < l_L < 4L$$
such that $[\bx_{l_1}]^{-1},\dots,[\bx_{l_L}]^{-1}$ is a basis for $Z$ as a right $D$-vector space, which satisfies \eqref{basis_1}; we will write $\bwy_n = [\bx_{l_n}]^{-1}$ for each $1 \leq n \leq L$. Notice that in fact $\bwy_1,\dots,\bwy_L \in \O_D^N$, where $\O_D$ is defined in~\eqref{O_D}. Combining Lemmas~3.4 and~3.5 of \cite{quaternion}, we see that
\begin{equation}
\label{ht1}
H([Z]) = H^{\O_D}(Z)^4 \leq \MM(\O)^{4(N-L)} H^{\O}(Z)^4,
\end{equation}
while Lemma~3.1 of \cite{quaternion} implies that for each $1 \leq n \leq 4L$
\begin{equation}
\label{ht2}
h([\bx_n]^{-1}) \leq 2s(\alpha,\beta) h(\bx_n).
\end{equation}
Combining \eqref{gen_bnd_1} with \eqref{ht1} and \eqref{ht2} yields
$$h(\bwy_n) \leq 4 L (M+2J+1)^{\frac{1}{d}} |\D_K|^{\frac{L+1}{2d}} s(\alpha,\beta) \MM(\O)^{4(N-L)} H^{\O}(Z)^4$$
for each $1 \leq n \leq L$. Arranging $\bwy_1,\dots,\bwy_L$ in the non-decreasing height order yields~\eqref{mn1} and completes the proof of Theorem~\ref{main1}.
\endproof

\begin{rem} Theorem~\ref{main1} is a version of Theorem~A.1 of \cite{quad_zero} over a quaternion algebra. It constitutes a non-commutative version of Sielgel's lemma missing a union of varieties and hence generalizes a non-commutative version of Siegel's lemma first established by Liebend\"orfer in \cite{liebendorf:1}.
\end{rem}

\begin{thm} \label{main2} Let all the notation be as in Theorem~\ref{main1}, and let $F(\bX,\bY) \in D[\bX,\bY]$ be a hermitian form in $2N$ variables. Suppose that there exists a point $\bwy \in Z \setminus \left( \left( \bigcup_{m=1}^M U_m \right) \left( \bigcup_{l=1}^J W_l \right) \right)$ such that $F(\bwy) := F(\bwy,\bwy) = 0$, then there exists such a point with
\begin{equation}
\label{mn2}
h(\bwy) \leq \A_{K,\O} (L,M,J,\alpha,\beta) H_{\inf}(F)^{\frac{9L+11}{2}} H^{\O}(Z)^{4(9L+12)},
\end{equation}
where the constant $\A_{K,\O} (L,M,J,\alpha,\beta)$ given by~\eqref{AK} above. Furthermore, there exists a point $\bz \in D^N \setminus \left( \bigcup_{m=1}^M U_m \right)$ such that $F(\bz)=0$ and
\begin{equation}
\label{mn3}
h(\bz) \ll_{K,N,M} 2s(\alpha,\beta) \left( \frac{2s(\alpha,\beta)^2}{t(\alpha,\beta)} H_{\inf}(F) \right)^{\frac{N+1}{2}}.
\end{equation}
\end{thm}

\proof
Notice that
$$[\bwy] \in [Z] \setminus \left( \left( \bigcup_{m=1}^M [U_m] \right) \left( \bigcup_{l=1}^J [W_l] \right) \right)$$
and $Q_F([\bwy]) = 0$. Let $\W$ be the Witt index and $\lambda$ the dimension of the radical of the quadratic space $([Z],Q_F)$ over $K$, so that a maximal totally isotropic subspace of $([Z],Q_F)$ has dimension $\mu := \W+\lambda$, then Theorem~1.1 of \cite{quad_zero} guarantees that there exist $\mu$ linearly independent vectors
$$\bx_1,\dots,\bx_{\mu} \in [Z] \setminus \left( \left( \bigcup_{m=1}^M [U_m] \right) \left( \bigcup_{l=1}^J [W_l] \right) \right)$$
such that for each $1 \leq n \leq \mu$
\begin{equation}
\label{ht3}
h(\bx_n) \leq T_K(L,M+2J+1) H(Q_F)^{\frac{9L+11}{2}} H([Z])^{9L+12},
\end{equation}
where $T_K(L,M)$ is a dimensional field constant, given by equation (43) of \cite{quad_zero}; its technical definition is somewhat complicated, so we do not present here in the interest of the brevity of exposition. Now, combining \eqref{ht3} with \eqref{ht1}, \eqref{ht2} and Lemma~3.2 of \cite{quaternion}, we obtain
\begin{equation}
\label{ht4}
h([\bx_n]^{-1}) \leq \A_{K,\O} (L,M,J,\alpha,\beta) H_{\inf}(F)^{\frac{9L+11}{2}} H^{\O}(Z)^{4(9L+12)}
\end{equation}
for each $1 \leq n \leq \mu$. Since
$$[\bx_1]^{-1},\dots,[\bx_{\mu}]^{-1} \in Z \setminus \left( \left( \bigcup_{m=1}^M U_m \right) \left( \bigcup_{l=1}^J W_l \right) \right),$$
and $\mu \geq 1$, we can take, for instance, $\bwy = [\bx_1]^{-1}$, and obtain \eqref{mn2}. 

Finally, to obtain \eqref{mn3}, we can combine Theorem of \cite{dietmann} with \eqref{ht2} and Lemma~3.1 of~\cite{quaternion} in the same manner as above. This completes the proof of Theorem~\ref{main2}.
\endproof

\begin{rem} \label{rem1} Inequality \eqref{mn2} of Theorem~\ref{main2} is a version of Theorem~1.1 of \cite{quad_zero} and \eqref{mn3} is a version of the main theorem of \cite{dietmann} (see also \cite{me:smallzeros}), both over a quaternion algebra. The bound of \eqref{mn3} demonstrates better dependence on $H_{\inf}(F)$ when  $Z=D^N$, although it only provides a point missing a collection of linear subspaces.  
\end{rem}

One can continue applying our ``transfer method" in the similar manner to obtain analogues of results on Siegel's lemma outside of linear subspaces with an additional dependence on the height of these subspaces (see \cite{me:number}, \cite{gaudron}).
\bigskip

\bibliographystyle{plain}  
\bibliography{fukshansky_henshaw}        

\end{document}